\documentstyle[12pt,amsfonts,amssymb,leqno]{article}

\newfont\got{eufm10}

\newtheorem{proposition}{Proposition}[section]

\newtheorem{lemma}[proposition]{Lemma}
\newtheorem{defn}[proposition]{Definition}

\newcounter{secnum}

\begin{document}
\setcounter{section}{-1}

\begin{center}
{\Large \bf Core models in the presence of Woodin
cardinals}
\end{center}

\begin{center}
{\Large \tt 
by Ralf Schindler}
\end{center}

\begin{center}
{\large \tt 
www.logic.univie.ac.at/${}^\sim$rds}
\end{center}

\section{Introduction.}

In this paper we show that the core model might exist even if
there are Woodin cardinals in $V$.
This observation is not new.
Woodin \cite{hugh}, in his proof that ${\sf AD}_{\mathbb R}$ 
implies the ${\sf AD}_{\mathbb R}$ hypothesis, 
constructed models of ${\sf ZFC}$ in
which there are fully iterable extender models with Woodin cardinals
which satisfy (among other things) 
a weak covering property. Steel \cite{john},
in his proof that $M_n$ satisfies
$V = HOD$, gave an argument which appears 
to be a special case of what we shall do in this paper.
However, the general method for
constructing the core model in the theory
${\sf ZFC}$ + ``there is a measurable cardinal above $n$ Woodin cardinals''
+ ``$M_{n+1}^\#$ does not exist'' 
which we shall present here
does seem to be new. 

This method might in turn admit generalizations, but we do not know how
to do it.
We shall indicate that our method might have applications;
we shall prove that if $M$ is an ultrapower of $V$ by an extender
with countably closed support then $K^M$ is an iterate of $K$ (although
$K$ might not be fully iterable). 

\begin{defn}
We let $\Sigma$ denote the following partial function.
$\Sigma({\cal T}) = b$ if and only if ${\cal T}$ is an iteration tree of
limit length and either

(a) ${\rm cf}({\rm lh}({\cal T})) > \omega$ and $b$ is the unique
cofinal branch through ${\cal T}$, or else

(b) 
$b$ is the unique cofinal branch through ${\cal T}$
such that ${\cal Q}(b,{\cal T})$ exists and is weakly iterable.

We say that an iteration tree ${\cal T}$ is {\em according to} $\Sigma$ if and
only if $[0,\lambda)_T = \Sigma({\cal T} \upharpoonright \lambda)$ for
every limit ordinal $\lambda < {\rm lh}({\cal T})$.
\end{defn}

\begin{defn}
Let ${\cal M}$ be a $k$-sound premouse, and let $\delta \in {\cal M}$.
We say that ${\cal M}$ is {\em weakly iterable above} $\delta$ if
for all weak $k$-embeddings $\pi \colon {\bar {\cal M}} \rightarrow 
{\cal M}$ with $\delta \in {\rm ran}(\pi)$,
$\Sigma$ witnesses that ${\bar {\cal M}}$ is $(k,\omega_1+1)$ iterable
above $\pi^{-1}(\delta)$ (i.e., 
with respect to iteration trees in which all extenders
used have critical point $\geq \pi^{-1}(\delta)$).
If $\delta = 0$ then we omit ``above $0$.''
\end{defn}

By arguments of
Woodin and Neeman, 
the $K$ of the next section can be shown to be not ``fully iterable.'' 
We have to introduce a new kind of iterability in order to be able to formulate to which
extent $K$ will be iterable after all.

\begin{defn}\label{delta-iterability}
Let ${\cal M}$ be a premouse, and let ${\vec \delta} = (\delta_1 < \delta_2 < 
... < \delta_n)$ be a sequence of inaccessible cardinals.
We say that ${\cal M}$ is ${\vec \delta}$ iterable via $\Sigma$
provided the 
following holds true.

Suppose that ${\cal T}$ is a putative iteration tree on ${\cal M}$ which 
is according to $\Sigma$.
Suppose also that we can write $${\cal T} = {\cal T}_0^{ \ \frown} 
{\cal T}_1^{ \ \frown} ... {}^\frown
{\cal T}_n$$ where the following assumptions are met.

(a) For each $k \leq n$, ${\cal T}_k$ is an iteration tree
of length $\theta_k$ (possibly $\theta_k = 1$, i.e., ${\cal T}_k$ is trivial),

(b) for each $k < n$, $\theta_k$ is a successor ordinal and 
${\cal T}_{k+1}$ is an iteration tree on ${\cal M}^{{\cal T}_k}_{\theta_k-1}$,
the last model of ${\cal T}_k$,

(c) for each $k \leq n$, if $k > 0$ then
${\rm crit}(E_\alpha^{{\cal T}_k}) >
\delta_k$ whenever $\alpha+1 < {\rm lh}({\cal T}_k)$
and if $k < n$ then ${\rm lh}(E_\alpha^{{\cal T}_k}) < \delta_{k+1}$
whenever $\alpha+1 < {\rm lh}({\cal T}_k)$,
and

(c) for each $k < n$, either $\theta_k < \delta_{k+1}$ or else
$\theta_k \in \{ \delta_{k+1}, \delta_{k+1}+1 \}$ and there
is an unbounded $A_k \subset \delta_{k+1}$ and 
a non-decreasing sequence $(\Omega_i \colon i \in A_k)$
of inaccessible cardinals below $\delta_{k+1}$ converging to $\delta_{k+1}$
such that
if $j \leq i \in A_k$ then $j \leq_{T_k} i$, $\pi_{0 j}^{{\cal T}_k} \
{\rm " } \ \Omega_j \subset \Omega_j$, 
and $\pi_{ji}^{{\cal T}_k} \upharpoonright 
\Omega_j = {\rm id}$.

Then either ${\cal T}$ has successor length and its last model is well-founded,
or else ${\cal T}$ has limit length and $\Sigma({\cal T})$ is well-defined. 
\end{defn}

Even if ${\cal M} \cap {\rm OR} < \delta_n$ then it makes sense 
to say that ${\cal M}$ is ${\vec \delta}$ iterable. 
If ${\cal M} \cap {\rm OR} \geq \delta_k$ and if ${\cal T}$ is as above
then it is easy to see that $\delta_1$, ..., $\delta_k$ are not moved
by the relevant embeddings in ${\cal T}$.
If ${\cal T}_k \upharpoonright \delta_{k+1}$ has length $\delta_{k+1}$
then there is a unique cofinal branch through 
${\cal T}_k \upharpoonright \delta_{k+1}$ by Definition \ref{delta-iterability}
(c), and hence $\Sigma({\cal T}_k \upharpoonright \delta_{k+1})$ 
is then well-defined.

\section{The construction.}

Fix $n < \omega$.
Throughout this section we shall assume that $\delta_1 < \delta_2 < ... < 
\delta_n$ are Woodin cardinals, that $\Omega > \delta_n$ is measurable, 
and that $M_{n+1}^\#$ does not exist.
We aim to construct $K$, 
``the core model
of height $\Omega$.''

We sometimes write $\delta_0 = 0$ and $\delta_{n+1} = \Omega$.
We shall recursively construct $K||\delta_k$, $k \leq n+1$, in
such a way that the following assumptions
will hold inductively.

\bigskip
{\bf A}${}_{1,k}. \ $ Suppose that $k>0$.
If $\kappa$ is measurable with $\delta_{k-1} < \kappa \leq \delta_k$ then
$K||\kappa$ is the core model over $K||\delta_{k-1}$ of height $\kappa$.

\bigskip
{\bf A}${}_{2,k}. \ $ Suppose that $k>1$.
Then $\delta_1$, ..., $\delta_{k-1}$ are Woodin cardinals
in $K||\delta_k$.

\bigskip
{\bf A}${}_{3,k}. \ $ $K||\delta_k$ is ${\vec \delta}$ iterable.

\bigskip
Now suppose that $K||\delta_k$ has already been contructed, where
$k<n+1$. Suppose also that 
{\bf A}${}_{1,k}$, {\bf A}${}_{2,k}$, and {\bf A}${}_{3,k}$ hold.
We aim to construct $K||\delta_{k+1}$ in such a way that 
{\bf A}${}_{1,k+1}$, {\bf A}${}_{2,k+1}$, and {\bf A}${}_{3,k+1}$ hold.

Let us first run the $K^c(K||\delta_k)$ construction inside $V_{\delta_{k+1}}$.
Let ${\cal N}_\xi$ and ${\cal M}_\xi = 
{\frak C}_\omega({\cal N})$ be the models 
of this construction. We would let the construction break down just
in case that we reached some $\xi < \delta_{k+1}$ 
with $\rho_\omega({\cal N}_\xi)
< \delta_k$. 
We'll show in a minute that this is not the case.

We shall prove the following statements under the assumption that
$\xi \leq \delta_{k+1}$ and ${\cal N}_\xi$ is defined.

\bigskip
{\bf A}${}_{1,k,\xi}. \ $ Suppose that $\xi>0$.
Then $\delta_1$, ..., $\delta_k$ are Woodin cardinals in ${\cal N}_\xi$.

\bigskip
{\sc Proof} of {\bf A}${}_{1,k,\xi}$. By {\bf A}${}_{2,k}$ it suffices
to show that $\delta_k$ is a Woodin cardinal in ${\cal N}_\xi$, where 
we assume that $k>0$.
Let $f \colon \delta_k \rightarrow \delta_k$, $f \in {\cal N}_\xi$,
and let $g \colon \delta_k \rightarrow \delta_k$ be defined by
$g(\xi) =$ the least $V$-inaccessible strictly 
above $f(\xi) \cup \delta_{k-1}$. As $\delta_k$ 
is Woodin in $V$, there is some $\kappa < \delta_k$
and some $V$-extender $E \in V_{\delta_k}$ such that if
$\pi \colon V \rightarrow_E M$ then $M$ is transitive,
$V_{\pi(g)(\kappa)} \subset M$, and $\pi(K||\delta_k) \cap V_{\pi(\kappa)} =
K||\delta_k \cap V_{\pi(\kappa)}$.
{\bf A}${}_{1,k}$ then gives, by standard arguments,
a total extender $F \in K||\delta_k$
such that if $\sigma \colon K||\delta_k \rightarrow_F K'$ then
$K'$ is transitive and
$K \models V_{\sigma(f)(\kappa)} \subset K'$. \hfill $\square$ 
({\bf A}${}_{1,k,\xi}$) 

\bigskip
{\bf A}${}_{2,k,\xi}. \ $ 
${\cal N}_\xi$ is weakly iterable above $\delta_k$.

\bigskip
{\sc Proof} of {\bf A}${}_{2,k,\xi}$. Cf.~\cite[\S \S 1, 2, and 9]{CMIP}. 
\hfill $\square$ ({\bf A}${}_{2,k,\xi}$)

\bigskip
{\bf A}${}_{3,k,\xi}. \ $  
Suppose that there is no $\delta > \delta_k$ which is 
definably Woodin in ${\cal N}_\xi$.
Then
${\cal N}_\xi$ is ${\vec \delta}$ iterable via $\Sigma$.

\bigskip
{\sc Proof} of {\bf A}${}_{3,k,\xi}$.
Suppose not.
Let $${\cal T} = {\cal T}_0^{ \ \frown} 
{\cal T}_1^{ \ \frown} ... {}^\frown {\cal T}_n$$ be a putative iteration
tree on ${\cal N}_\xi$ according to $\Sigma$ and as in Definition
\ref{delta-iterability}
such that either ${\rm lh}({\cal T})$ is a 
successor ordinal and ${\cal T}$ has a last ill-founded model
or else ${\rm lh}({\cal T})$ is a limit ordinal and
$\Sigma({\cal T})$ is undefined. (Notice that 
if $k<n$ then
we must have that ${\cal T}_{k+1}$, ..., 
${\cal T}_n$ are trivial.)

Suppose for a second that ${\rm lh}({\cal T})$ is a limit ordinal.
Then ${\rm lh}({\cal T}) \notin \{ \delta_1, ..., \delta_n \}$, as otherwise
$\Sigma({\cal T})$ would clearly be defined by Definition 
\ref{delta-iterability} (c).
Set $\delta = \delta({\cal T})$,
and let $W = K^c({\cal M}({\cal T}))$ (built of
height ${\rm OR}$).
$W$ is well-defined and in fact
$\delta \notin \{ \delta_1 , ..., \delta_n \}$ 
is a Woodin cardinal in $W$ due to the fact that 
$\Sigma({\cal T})$ is undefined (cf.~\cite{CMMW}). 
By {\bf A}${}_{1,k,\xi}$ and by 
the proof of \cite[Theorem 11.3]{FSIT}, we'll also
have that $\delta_1$, ..., $\delta_n$ are Woodin cardinals
in $W$, and $\Omega$ will be
measurable in $W$. As $M_{n+1}^\#$ does not exist 
we know that $W||\Omega^{++}$ can't be weakly iterable.

If ${\rm lh}({\cal T})$ is a successor ordinal then we
let $\delta$ and $W$ be undefined.

Now let $\Theta$ be regular and large enough, and pick $${\bar H} 
\stackrel{\pi}{\rightarrow} H' \stackrel{\sigma}{\rightarrow}
H_\Theta$$ such that ${\bar H}$ and $H'$ are transitive, 
${\rm Card}({\bar H}) = \aleph_0$, and 
$\{ {\cal N}_\xi , {\cal T} \} \subset {\rm ran}(\sigma \circ \pi)$.
If $k>0$ then we shall also assume that 
${\rm Card}(H') < \delta_k$ and 
${\rm ran}(\sigma) \cap \delta_k \in \delta_k$. (If $k=0$ then 
we allow $\sigma = {\rm id}$.) 
Set ${\bar {\cal N}} = (\sigma \circ \pi)^{-1}({\cal N}_\xi)$ and
${\bar {\cal T}} = (\sigma^{-1} \circ \pi)({\cal T})$.
If ${\rm lh}({\cal T})$ is a limit ordinal then we shall
also assume that $W||\Omega^{++} \in {\rm ran}(\sigma \circ \pi)$, and
setting
${\bar W} = (\sigma \circ \pi)^{-1}(W||\lambda)$ we may and shall assume that
${\bar W}$ is not $\omega_1+1$ iterable.
We also write ${\cal N}' = \pi({\bar {\cal N}}) = \sigma^{-1}({\cal N}_\xi)$
and $\delta' = \sigma^{-1}(\delta_k)$ ($= {\rm crit}(\sigma)$ if 
$k>0$).

\bigskip
{\bf Claim 1.}
Suppose that $k>0$, and
let $\kappa < \delta_k$ be measurable with ${\cal N}' \cap {\rm OR}
< \kappa$.
Then ${\cal N}'$ is $\Sigma$ coiterable with $K||\kappa$.

\bigskip
{\sc Proof} of Claim 1. Suppose not. Let ${\cal W}$ and ${\cal W}'$ denote
the putative iteration trees arising from the coiteration
of ${\cal N}'$ with $K||\kappa$. ${\cal W}$ is above $\delta'$ and
${\cal W}'$ is above $\delta_{k-1}$. We assume
that ${\cal W}$ has a last ill-founded model or else that
$\Sigma({\cal W})$ is not defined. By {\bf A}${}_{2,k,\xi}$, we must have
that ${\cal W}$ has limit length and $\Sigma({\cal W})$ is not defined.
By {\bf A}${}_{3,k}$, $\Sigma({\cal W}')$ is well-defined, i.e.,
${\cal Q}({\cal W}')$ does exists.  

%
Let us pick $$\tau \ \colon \ 
H'' \rightarrow H_\Theta$$ such that $H''$ is transitive,
${\rm Card}(H'') = \aleph_0$, and $\{ {\cal N}', {\cal W}, 
{\cal Q}({\cal W}')
\} \subset {\rm ran}(\tau)$.
Set ${\bar {\cal N}}' = \tau^{-1}({\cal N}')$ and
${\bar {\cal W}} = \tau^{-1}({\cal W})$.

As ${\cal N}_\xi$ is weakly iterable above $\delta_k$,
we certainly have that ${\bar {\cal N}}'$ is $\omega_1+1$ iterable
above $\tau^{-1}(\delta')$.
In particular, ${\cal Q}({\bar {\cal W}})$ exists, and
${\cal Q}({\bar {\cal W}}) 
\trianglelefteq {\cal M}_{\Sigma({\bar {\cal W}})}^{\bar {\cal W}}$. 
However, 
we must have that ${\cal Q}({\bar {\cal W}}) = 
\tau^{-1}({\cal Q}({\bar {\cal W}'}))$. In particular, 
${\cal Q}({\bar {\cal W}}) \in H''$, and
therefore by standard arguments
$\Sigma({\cal W})$ is well-defined after all. \hfill $\square$ (Claim 1)

\bigskip
{\bf Claim 2.} ${\bar {\cal N}}$ is $\omega_1+1$ iterable.

\bigskip
{\sc Proof} of Claim 2. 
Well, if $k=0$ then this readily follows from 
{\bf A}${}_{2,k,\xi}$.  
Let us now assume that $k>0$.
Let $\kappa < \delta_k$ be measurable with ${\cal N}' \cap {\rm OR}
< \kappa$.
By Claim 1,
we may let ${\cal W}$ and ${\cal W}'$ denote the
iteration trees arising from the (successful)
$\Sigma$ coiteration of ${\cal N}'$ with $K||\kappa$.
By {\bf A}${}_{1,k}$, we shall have that
$\pi_{0 \infty}^{{\cal W}} \colon 
{\cal N}' \rightarrow {\cal M}^{{\cal W}'}_\infty ||\gamma$ for some $\gamma$,
where $\pi_{0 \infty}^{{\cal W}}
\upharpoonright \delta' = {\rm id}$.
We then have that $$\pi_{0 \infty}^{{\cal W}} \ \circ \ \pi \ \colon \
{\bar {\cal N}} \rightarrow {\cal M}_\infty^{{\cal W}'} || \gamma$$
witnesses, by {\bf A}${}_{3,k}$, that
${\bar {\cal N}}$ is $\omega_1+1$ iterable via $\Sigma$.
\hfill $\square$ (Claim 2)

\bigskip
In particular, 
${\cal T}$ must have
limit length.
By Claim 2,
${\cal Q}({\bar {\cal T}})$ exists (and is $\omega_1+1$ iterable). Set
${\cal Q} = {\cal Q}({\bar {\cal T}}) 
\trianglelefteq {\cal M}_{\Sigma({\bar {\cal T}})}^{{\bar {\cal T}}}$, and
set ${\bar \delta} =
(\sigma \circ \pi)^{-1}(\delta)$. We have that ${\bar \delta}$ 
is a cutpoint in ${\bar W}$,
${\bar \delta}$ is a Woodin cardinal in ${\bar W}$, and
${\bar W}$ is $\omega_1+1$ iterable above ${\bar \delta}$.
Of course, ${\cal Q} || {\bar \delta} = {\bar W}||{\bar \delta}$.
Let ${\cal V}'$ and ${\cal V}$ denote the iteration trees arising from the
(successful) 
comparison of ${\bar {\cal Q}}$ with ${\bar W}$.
We shall have that $\pi_{0 \infty}^{\cal V} \colon {\bar W} 
\rightarrow {\cal M}^{{\cal V}'}_\infty || \gamma'$ for some $\gamma'$.
By Claim 2, this embedding
witnesses that
${\bar W}$ is $\omega_1+1$ iterable.
Contradiction!
\hfill $\square$ ({\bf A}${}_{3,k,\xi}$)
 
\bigskip
{\bf A}${}_{4,k,\xi}. \ $ 
$\rho_\omega({\cal N}_\xi) \geq \delta_k$.

\bigskip
{\sc Proof} of {\bf A}${}_{3,k,\xi}$. We may trivially
assume that $k>0$. Suppose that $\rho =
\rho_\omega({\cal N}_\xi)
< \delta_k$. Let $\Theta$ be large enough,
and pick $\pi \colon {\bar H} \rightarrow H_\Theta$ such that
${\bar H}$ is transitive, ${\rm Card}({\bar H}) < \delta_k$,
${\rm ran}(\pi) \cap \delta_k \in (\delta_k \setminus (\rho+1))$, and
$\{ {\cal N}_\xi , \delta_k \} \subset {\rm ran}(\pi)$.
Set ${\bar {\cal N}} = \pi^{-1}({\cal N}_\xi)$
and ${\bar \delta} = \pi^{-1}(\delta_k)$.
Let $\kappa < \delta_k$ be measurable such that $\kappa > {\bar {\cal N}}
\cap {\rm OR}$.
By {\bf A}${}_{3,k,\xi}$ (or rather by the proof of
Claim 1 in the proof of {\bf A}${}_{3,k,\xi}$), 
we may now successfully coiterate ${\bar {\cal N}}$ with $K||\kappa$.
{\bf A}${}_{1,k}$ implies that any witness to $\rho = \rho_\omega({\bar {\cal N}})$ is an element of $K||\kappa$. But any such witness
is also a witness to
$\rho = \rho_\omega({\cal N}_\xi)$. Contradiction! 

\hfill $\square$ 
({\bf A}${}_{3,k,\xi}$)

\bigskip
We have shown that ${\cal N}_{\delta_{k+1}}$ exists, is ${\vec \delta}$
iterable, and that $\delta_1$, ..., $\delta_k$ are Woodin cardinals
in ${\cal N}_{\delta_{k+1}}$. Let us write $K^c(K||\delta_k)$ for
${\cal N}_{\delta_{k+1}}$.

Now let $\Omega$ be a measurable
cardinals with
$\delta_k < \Omega \delta_{k+1}$.
Using $K^c(K||\delta_k)||\Omega$
we may define $K_\Omega$ 
as the core model over $K||\delta_k$. We shall have the
canonical embedding $\pi_\Omega \colon K_\Omega \rightarrow
K^c(K||\delta_k)||\Omega$ (cf.~\cite{maximality}).
If
$\Omega$ and $\Omega'$ are measurable
cardinals with
$\delta_k < \Omega \leq \Omega' \leq \delta_{k+1}$
then $K_\Omega \trianglelefteq K_{\Omega'}$ and 
$\pi_\Omega = \pi_{\Omega'} \upharpoonright K_\Omega$.
Let us define $K||\delta_{k+1}$ as the ``union''
of the $K_\Omega$, where $\Omega$ is a measurable cardinal
with $\delta_k < \Omega \delta_{k+1}$.
The ``union'' of the maps $\pi_\Omega$ then give a canonical
embedding $\pi \colon K||\delta_{k+1} \rightarrow K^c(K||\delta_k)$.
It is now easy to see that
{\bf A}${}_{1,k+1}$, {\bf A}${}_{2,k+1}$, and {\bf A}${}_{3,k+1}$ hold true.

\section{An application.}

\begin{lemma}
Suppose that $\delta_1 < ... < \delta_n$ are Woodin cardinals,
that $\Omega > \delta_n$ is measurable, and that $M_{n+1}^\#$ does
not exists. Let $E \in V_\Omega$ be an extender with countably 
closed support. Let $\pi \colon V \rightarrow_E M$ where $M$ is transitive.

Let $K$ denote the core model of height $\Omega$ constructed in
the previous section. There is then an iteration tree
${\cal T}$ on $K$ as in Definition \ref{delta-iterability} such that
${\cal M}_\infty^{\cal T} = \pi(K)$ and $\pi_{0 \infty}^{\cal T} =
\pi \upharpoonright K$. 
\end{lemma} 

{\sc Proof.} This follows by the methods of \cite{iterates-of-K}. \hfill
$\square$

\end{document}